\documentclass[12pt]{article}
\usepackage{amssymb,amsmath,amsthm,amscd,latexsym}
\usepackage{mathrsfs}
\usepackage{mathrsfs}
\usepackage{amsfonts}
\usepackage{amsmath}
\usepackage{amssymb}
\usepackage{amscd}
\usepackage{mathrsfs}
\usepackage{amssymb}
\usepackage{amsmath}
\usepackage{amsthm}
\usepackage{color}
\usepackage{latexsym}
\usepackage{indentfirst}
\usepackage{anysize}
\usepackage[normalem]{ulem}

\renewcommand{\paragraph}{\roman{paragraph}}
\input cyracc.def

\newcommand{\F}{\mathbb{F}}

\newtheorem{theorem}{Theorem}

\newtheorem{lem}{Lemma}
\newtheorem{proposition}[theorem]{Proposition}
\newtheorem{thm}{Theorem}

\newtheorem{conj}{Conjecture}

\newtheorem{coro}{Corollary}
\newtheorem{ex}{Example}

\theoremstyle{definition}
\newtheorem{definition}[theorem]{Definition}

\newcommand{\R}{\mathbb{R}}

\newcommand{\Z}{\mathbb{Z}}

\begin{document}
%\begin{CJK*}{GBK}{song}\CJKtilde
\title{\bf Permutations designs, projective planes, and transitive sets of permutations
\thanks{This research is supported by National Natural Science Foundation of China (12071001, 61672036), Excellent Youth Foundation of Natural Science Foundation of Anhui Province (1808085J20), the Academic Fund for Outstanding Talents in Universities (gxbjZD03).
}}

\author{
\small{Minjia Shi$^{1,3}$, {XiaoXiao Li$^{1}$}, Patrick Sol\'e$^{2}$}\\ %Patrick Sol\'e$^4$
\and %\small{${}^1$School of Mathematical Sciences, Anhui University, Hefei, 230601, China}\\
\small{${}^1$Key Laboratory of Intelligent Computing \& Signal Processing of Ministry of Education,}\\
\small{School of Mathematical Sciences, Anhui University,}\\
\small{Hefei 230601,  China}\\
  %\and \small{${}^2$Telecom ParisTech, Palaiseau, France}
\and
\small{${}^2$I2M,(Aix-Marseille Univ., Centrale Marseille, CNRS), Marseilles, France}\\
\small{${}^3$ Corresponding Author}
}

%Key Laboratory of Intelligent Computing and Signal Processing of Ministry of Education, School of Mathematics Sciences, Anhui University, Hefei, 230601, China.

\date{}
\maketitle
\begin{abstract} A notion of $t$-designs in the symmetric group on $n$ letters, called permutation designs in this paper, was introduced by Godsil in 1988. 
In particular $t$-transitive sets of permutations form a $t$-design in that sense. 
We derive special lower bounds
for $t=1$ and $t=2$ by a power moment method. For general $n,t$  we give a 
%linear programming lower bound . For  $n\ge 4$ and $t=2,$ this bound is strong enough to show 
a lower bound on the size of such $t$-designs of  $n(n-1)\dots (n-t+1),$ which is best possible when sharply $t$-transitive sets of permutations exist.
This shows, in particular, that tight $2$-designs do not exist.
\end{abstract}

{\bf Keywords:}  Permutations, Charlier polynomials, Latin squares, Projective planes\\

{\bf AMS Math Sc. Cl. (2010):} Primary 05E30, Secondary 05B15
%%%%%%%%%%%%%%%%%%%%%%%%%%%%%%%%%%%%%%%%%%%%%%%%%%%%%%%%%
\section{Introduction}
%%%%%%%%%%%%%%%%%%%%%%%%%%%%%%%%%%%%%%%%%%%%%%%%%%%%%%
In his celebrated thesis \cite{D}, Delsarte developed a theory of designs in $Q$-polynomial association schemes \cite{BCN}, and a linear programming bound on their sizes. 
This theory was extended in the 1980s by
Godsil with the notion of designs in polynomial spaces \cite{G,Gp}. This latter theory, in particular, introduced a notion of permutation designs, or designs in the symmetric group,
which was developed further in \cite{CG}.
In particular, $t$-transitive groups and $t$-transitive sets of permutations are $t$-designs in that sense. The existence of permutation designs not coming from transitive sets of 
permutations is a challenging open problem. This
notion of design cannot be understood in the context of $Q$-polynomial association schemes, as the conjugacy scheme of the symmetric group is not $Q$-polynomial.
Unfortunately, the theory of Polynomial Spaces does not comprise  a concept of inner distribution or dual inner distribution as $Q$-polynomial association schemes do. It does contain a linear 
programming bound \cite[ Theorem 14.5.3]{G}, but this bound is written in terms of polynomials. 
 In a recent paper \cite{SRS}, a theory of
Distance Degree Regular metric spaces (DDR spaces) was introduced, with a notion
of frequencies playing the role of the inner distribution in classical $Q$-polynomial association schemes, and a concomitant notion of dual frequencies based on a system of orthogonal polynomials 
attached to the DDR space under consideration. This allows for a criterion for a set of permutations to form a design (Proposition \ref{crit}). While   the 
relevant orthogonal polynomials are not known in general,
in the special case of the symmetric group it is
possible to use the work of Tarnanen on the conjugacy scheme of the symmetric group \cite{T} to define these polynomials in degree at most $n/2.$
In that remarkable work, Tarnanen provides an orthogonality relation for Charlier polynomials that cannot be found in the classical treatises
on orthogonal polynomials \cite{KS,Sz}. We make explicit the dual frequency criterion when $t=1,2.$ For general $t,$ we derive a lower bound of $n(n-1)\dots (n-t+1),$ 
on the size of a $t$-design (Theorem \ref{sm}).
This is best possible when sharply $t$-transitive sets of permutations exist. In particular, such sets are known to exist for $t=2$ when a projective plane of order $n$ exists.
Note that Theorem \ref{sm} is stronger than Corollary \ref{cor2}, which implies that tight $2$-designs in the polynomial space of $S_n$ do not exist.

The material is arranged as follows. The next section collects the necessary notations and definitions on permutations and DDR spaces. Section 3 collects some basic results on DDR spaces.
Section 4 studies in particular $1$-designs and $2$-designs by a power moment method. %Section 5 contains the linear programming bound and its numerical examples.
%Sections $3$ and $4$ recall some results established in \cite{SRS}. Section 5 establishes the linear 
%programming bound for designs in DDR spaces, and deduces from that the definition and existence condition  of tight designs. Section 6 specializes these results to the DDR space of permutations.
Section 5 concludes the paper.
%%%%%%%%%%%%%%%%%%%%%%%%%%%%%%%%%%%%%%%%%%%%%%%%%%%%
\section{Background material}
\subsection{Permutations and designs}

A permutation group $G$ acting on a set $X$ of $n$ elements is transitive if there is only one orbit on $X.$ It is $t$-transitive if it is transitive in its action on ${ X \choose t}$ the set of distinct $t$-uples from $X.$
It is sharply $t$-transitive if this action is regular, concretely if $|G|=\frac{n!}{(n-t)!}.$ We extend this terminology by relaxing the group hypothesis to a set of permutations acting on $X.$
It is well-known amongst geometers and group theorists  that a set of sharply $2$-transitive permutations on a set of size $n$
is equivalent to the existence of a projective plane $PG(2,n),$ that is to say a
$2-(n^2+n+1,n+1,1)$ design \cite{C}.

\subsection{DDR spaces}

For any pair of nonegative integers $\mathcal{M},\mathcal{N},$ denote by  $[\mathcal{M}..\mathcal{N}]$  the set of integers in the range $[\mathcal{M},\mathcal{N}].$
A finite metric space $(X,d)$ is {\em distance degree regular} (DDR) if its distance degree sequence is the same for every point.
Assume $(X,d)$ to be of {\em diameter} $n.$ 
In that case $(X,d)$ is DDR iff for each $0\le i\le n$ the graph $\Gamma_i=(X,E_i)$ which connects vertices at distance $i$ in $(X,d)$ is regular of {\em degree } $v_i.$ 
Thus $E_0=\{(x,x) \mid x \in X\}$ is the  diagonal of $X^2.$
Note that the $E_i$'s form a partition of $X^2.$
Examples of DDR spaces are as follows. In the first two examples the metric is the shortest path distance defined by a graph on $X.$
We assume that the reader has some familiarity with the theory of distance-regular (DR) graphs as can be found, for instance, in \cite{BI,BCD,BCN}.
Every DR graph on a set X induces a DDR metric space by using the shortest path distance of that graph.
\begin{enumerate}
\item The {\em Hamming graph} $H(n,q)$ is a graph on $\F_q^n$ two vertices being connected if they differ in exactly one coordinate. This graph is DR with valencies
$$v_i={n \choose i} (q-1)^i.$$ For examples of application of this graph to Coding Theory see \cite[chap. 17 \S 7,chap. 21]{MS}.
\item The {\em Johnson graph} $J(\nu,d)$  is a graph on the subsets of cardinality $d$ of a set of cardinality $\nu.$ (Assume $2d<\nu$).
Two subsets are connected iff they intersect in exactly $d-1$ elements.
This graph is DR with valencies
$$v_i={d \choose i}{\nu-d \choose i} .$$
Note that $J(\nu,d)$ can be embedded in $H(\nu,2)$ by identifying subsets and characteristic vectors.
For applications of this graph to design theory ( e.g. tight designs in the classical sense) see \cite{B,BI,DS}.
\item Consider the {\em symmetric group} on $n$ letters $S_n$ with metric
$$d_S(\sigma,\theta)=n-F(\sigma \theta^{-1}),$$
where $F(\nu)$ denotes the number of fixed points of $\nu.$
The space $(S_n,d_S)$ is a DDR metric space.
Let $w_k$ denote the numbers of permutations on $n$ letters with $k$ fixed points. A generating function for these numbers (sometimes called rencontres numbers) is
$$\sum_{k=0}^nw_ku^k=n!\sum_{j=0}^n\frac{(u-1)^j}{j!},$$ as per \cite{w}. Note that $v_i=w_{n-i}.$ It is clear that $d_S$ is not a shortest path distance since $ d_S(\sigma,\theta)=1$ is impossible.
Codes in $(S_n,d_S)$ were studied in \cite{T} by using the conjugacy scheme of the group $S_n.$ 
%However, in contrast with the next two subsections, this scheme is neither induced by a graph nor $Q$-polynomial.

\end{enumerate}
%%%%%%%%%%%%%%%%%%%%%
\section{Preliminaries}
%%%%%%%%%%%%%%%%%%%%%%%%%%%%%%%%%%%%%%%%%5
If $D$ is any non void subset of $X$ we define its {\em frequencies} as
$$\forall i \in [0..n], \, f_i=\frac{|D^2\cap E_i|}{|D|^2}.$$ Thus $f_0=\frac{1}{|D|},$ and $\sum\limits_{i=0}^n f_i=1.$ Note also that if $D=X,$ then $f_i=\frac{v_i}{|X|}.$
%Consider the random variable $a_D$ defined on $D^2$ with values in $[0..n]$ which affects to an equiprobably chosen 
%$(x,y)\in D^2$ the only $i$ such that $(x,y)\in E_i.$ Thus the frequencies $f_i=Prob(a_D=i).$ 

{\definition The set $D \subseteq X$ is a {\em $t$-design} for some integer $t$ if

$$\sum_{j=0}^n f_jj^i=          \sum_{j=0}^n \frac{v_j}{v}j^i.           $$
%$$\EE(a_D^i)=\EE(a_X^i)$$ 
for $i=1,\dots,t.$
}
%Let $z_0,z_1,\cdots,z_d$ denote arbitrary pairwise distinct real numbers.

(Note that trivially $\sum\limits_{j=0}^n f_jj^0=1$ so that we do not consider $i=0.$)
Thus, distances in $t$-designs are very regularly distributed. For a $2$-design, for instance, the average and variance of the distance coincide with that of the whole space.
 We will see in the next section that in the case of Hamming and Johnson graphs,
we obtain  classical combinatorial objects: block designs, orthogonal arrays. 

{\bf Remark:} In the three examples of Hamming, Johnson graphs, and the symmetric group our notion of design coincides with that of design in the respective polynomial space of \cite{G}.

{\definition We define a scalar product on $\R[x]$ attached to $D$  by the relation
$$\langle f,g \rangle_D=\sum_{i=0}^n f_i f(i)g(i). $$
Thus, in the special case of $D=X$ we have
$$\langle f,g \rangle_X=\frac{1}{|X|}\sum_{i=0}^n v_i f(i)g(i). $$}
We shall say that a sequence  $\Phi_i(x)$  of  polynomials of degree $i$ is {\em orthonormal of size $N+1$} if it satisfies
$$\forall i,j \in [0..N],\,  \langle \Phi_i,\Phi_j \rangle_X =\delta_{ij},$$
where $N\le n, $ the letter $\delta$ denotes the Kronecker symbol. That sequence is uniquely defined if we assume the leading coefficient of all $\Phi_i(x)$ for $i=0, 1,\dots, N$ to be positive.

%The following Proposition is given without proof in \cite[p.26]{Sz}.

%{\proposition \label{fonda} If at least $N+1$ of the $v_i$'s are nonzero, then $\langle , \rangle_X$ admits an orthormal system of polynomials of degree $N.$
For a given DDR metric space $(X,d),$ we shall denote by $N(X)$ the largest possible such $N.$ For instance if $X$ is an $n$-class $P$- and $Q$-polynomial association scheme, it is well-known that $N(X)=n.$
This fact is extended to DDR graphs in the next Proposition from \cite{SRS}.

{\proposition (\cite{SRS} )\label{fonda} If none of the $v_i$'s are zero, then $\langle , \rangle_X$ admits an orthonormal system of polynomials of size $n+1.$
In particular, the metric space of a DDR graph admits an orthonormal system of polynomials of size $n+1.$

%For a given DDR metric space $(X,d),$ we shall denote by $N(X)$ the largest possible such $N.$ For instance if $X$ is an $n$-class $Q$-polynomial associations scheme, it is well-known that $N(X)=n.$
}

{\definition For a given $D\subseteq X$ the {\em dual frequencies} are defined for $i=0,1,\dots,N(X)$ as $$\widehat{f_i}=\sum_{k=0}^n\Phi_i(k)f_k.$$}

%{\definition For a given $D\subseteq X$ the {\em cumulative distribution function} (c.d.f.) is defined as $$F_D(x)=Prob(a_D\le x)=\sum\limits_{i\le x}f_i.$$}

%%%%%%%%%%%%%%%%%%%%%%%%%%%%%%%%%%%%%%%%%%%%%%%%%%%%%%%%%
\section{Structure theorems}
First, we recall the characterization of $t$-designs in terms of dual frequencies.
{\proposition \label{crit} (\cite{SRS}) Let $t$ be an integer $\in [1..N(X)].$ The set $D\subseteq X$ is a $t$-design iff $\widehat{f_i}=0$ for $i=1,\dots,t.$}

To make this criterion more concrete we give an explicit expression for the $\Phi_i$ in terms of Charlier polynomials.

 Let $$C_k(x)=(-1)^k+\sum\limits_{i=1}^k(-1)^{k-i}{k \choose i} x(x-1)\cdots(x-i+1).$$ An exponential generating function is
given in \cite[(1.12.11)]{KS} as:
$$e^t(1-t)^x=\sum_{n=0}^\infty C_n(x)\frac{t^n}{n!}. $$  

Thus, for concreteness, $C_0(x)=1,\,C_1(x)=x-1,\, C_2(x)=x^2-3x+1.$

%Denote by $F(\sigma, \tau)$ the number of fixed points of $\sigma\tau^{-1}.$ let $d=n-F$ denote the Hamming distance on the permutation images.
The scalar product attached to the DDR space $(S_n,d_S)$ is then 
 $$\langle f,g\rangle_n=\frac{1}{n!} \sum_{k=0}^n w_{n-k} f(k)g(k).$$
 
 It is remarkable that the following orthogonality relation is not found in the classical treatises \cite{KS,Sz} on orthogonal polynomials.

{\lem \label{Charlier} The reversed Charlier polynomials $\widehat{C_k(x)}=C_k(n-x)$ satisfy the orthogonality relation $$ \langle\widehat{C_r},\widehat{C_s}\rangle_n=r! \delta_{rs},$$
for $r,s\le n/2.$}

\begin{proof}
By  Corollary 1 of \cite{T} we know that the $C_k(x)$'s satisfy the orthogonality relation 
$ (C_r,C_s) =r! \delta_{rs},$ for $r,s\le n/2.$
w.r.t. the inner product

 $$(f,g)_n=\frac{1}{n!} \sum_{k=0}^n w_k f(k)g(k).$$
 
 By the change of variable $x\mapsto n-x$ in $C_k(x)$ the result follows.
 
\end{proof}

%$$\sum_{k=0}^n f_kk^i=\sum_{k=0}^n \frac{v_k}{v}k^i,$$ for $i=1,2\dots,t,$ which, in turn, is equivalent to
%$$\EE(a_D^i)=\sum_{k=0}^n f_kk^i=\sum_{k=0}^n \frac{v_k}{v}k^i=\EE(a_X^i),$$ for $i=1,2\dots,t.$
%Hence, by linearity of the two scalar products, $\EE(a_D^i)=\EE(a_X^i)$ for $i=1,2\dots,t.$

Next, we connect the notion of designs in $Q$-polynomial association schemes with our notion of designs in metric spaces.
{\theorem (\cite{SRS}) If $(X,d)$ is the metric space induced by a $Q$-polynomial DR graph $\Gamma,$ with $z_k=k$ for $k=0,1,\dots,n,$ then 
a $t$-design in $(X,d)$ is exactly a $t$-design in the underlying association scheme of $\Gamma.$}

{\bf Examples:} The following two examples of interpretation of $t$-designs as classical combinatorial objects were observed first in \cite{D} and can be read about in \cite[chap. 21]{MS}.
\begin{enumerate}
\item If $\Gamma$ is the Hamming graph $H(n,q)$ then a $t$-design is an orthogonal array of strength $t.$ That means that every row 
induced by a  $t$-uple columns of $D$ sees the $q^t$ possible values a
constant number of times.
\item If $\Gamma$ is the Johnson graph $J(\nu,n)$ then a $t$-design $D$ is a combinatorial design of strength $t.$ This means the following. Consider $D$ as a collection of subsets of size $n$, traditionally called blocks. That means that every $t$-uple of elements of the groundset is contained in the same number $\eta$ of blocks.
    One says that $D$ is a $t-(\nu,n,\eta)$ design.
\end{enumerate}

Now, we give examples of $t$-designs in a metric space that is not a DR graph, or even a DDR graph. First, we give general links between $t$-designs and transitive sets of permutations.

{\theorem If $D\subseteq S_n$ is a $t$-transitive permutation group then it is a $t$-design in $(S_n,d_S).$ 
If $D\subseteq S_n$ is a $t$-design that is a subgroup of $S_n,$ then it is a $t$-transitive permutation group.
}

\begin{proof}
The first assertion is proved in \cite[Theorem 8]{SRS}. The second assertion follows by the method used in \cite[Th. 4.1]{CG} in the language of polynomial spaces. Let $s$ be a nonzero
integer $\le t.$
First, by Burnside Lemma the number of orbits of $Y$ on $[1..n]$ is $\langle 1, P\rangle _Y,$ with $P$ some polynomial of degree $s.$ Since $Y$ is a $t$-design we know that 
$\langle 1, P\rangle _Y=\langle 1, P\rangle_X.$ Since $S_n$ is $s$-transitive on $[1..n]$ we have that $\langle 1, P\rangle_X=1.$ The second assertion follows.
\end{proof}

Next, we give special results and constructions for $t=1,2.$

{\theorem \label{oned} A subset $D\subseteq S_n$ is a $1$-design in $(S_n,d_S)$  iff $\sum\limits_{j=0}^njf_j=n-1.$ In particular,
this condition is satisfied if we have $n$ permutations at pairwise distance $n$ when $f_1=f_2=\dots=f_{n-1}=0,$ and $f_n=\frac{n-1}{n}.$}

\begin{proof}
 The criterion can be obtained directly from the definition of a $1$-design by considering the generating function for the $w_k$'s given above.
 Alternatively,  if $n\ge 2,$ it can be derived by Proposition \ref{crit} above, and the expression for the first Charlier polynomial.
 The second assertion follows then by a direct calculation.
\end{proof}

 Theorem \ref{oned} yields a lower bound on the size of $1$-designs.

{\coro If $Y$ is a $1$-design in $(S_n,d_S)$ then $|Y|\ge n.$}

\begin{proof}
 Note that $$(n-1)=\sum\limits_{j=1}^njf_j\le n \sum\limits_{j=1}^nf_j=n(1-f_0)=n-\frac{n}{|Y|}.$$
\end{proof}

The existence of $n$ permutations of $S_n$ at pairwise Hamming distance $n$  is trivially equivalent to the existence of a {\bf Latin square} of order $n.$
This is the case when $Y$ is the group generated by a cycle of length $n.$ The Latin square is then the addition table of $(\Z_n,+).$\\

Here is a non-group example when $n=5,$ obtained from the smallest Latin Square that is not the multiplication table of a group.

$$Y=\{12345,24153,35421,41532, 53214\},$$ when $ 24153 \circ 35421 =13542 \notin Y.$

{\theorem \label{2crit} A subset $D\subseteq S_n$ is a $2$-design in $(S_n,d_S)$  iff $$\sum_{j=0}^njf_j=n-1, \& \sum_{j=0}^nj^2f_j=1+(n-1)^2.$$ In particular,
this condition is satisfied if we have $n(n-1)$ permutations with frequencies $f_1=f_2=\dots=f_{n-2}=0,$ and $ f_{n-1}=\frac{n-2}{(n-1)},\,f_{n}=\frac{1}{n}.$}

\begin{proof}
 The criterion can be obtained directly from the definition of a $2$-design by considering the generating function for the $w_k$'s given above.
 Alternatively, if $n\ge 4,$ it can be derived by Proposition \ref{crit} above, and the expression for the first and second Charlier polynomials.
 The second assertion follows then from the two equations of the criterion  by solving a 2 by 2 system for $ f_{n-1}$ and $f_n.$
\end{proof}

{\coro \label{cor2} If $Y$ is a $2$-design in $(S_n,d_S)$ then $|Y|\ge (n-1)^2+1.$}

\begin{proof}
{\bf First proof: elementary}
 Note that $$(n-1)=\sum\limits_{j=1}^njf_j=\sum\limits_{j=1}^n\sqrt{f_j}\sqrt{f_j}j\le \sqrt{1-f_0}\sqrt{1+(n-1)^2},$$
 where the inequality is a consequence of Cauchy-Schwarz inequality. Squaring both sides yields
 $$(n-1)^2 \le (1-f_0)(1+(n-1)^2)$$ hence $f_0(1+(n-1)^2)\le 1.$ The result follows by $f_0=1/|Y|.$\\
 
 {\bf Second proof: polynomial spaces} As observed in \cite[\S 4]{CG} $\dim Pol(\Omega,1)=1+(n-1)^2.$ The bound follows then by the Fisher-like inequality of  \cite[Th. 2.2]{CG}.
\end{proof}

Consider, as example of $2$-transitive group the affine group over $\F_n$ for $n$ a prime power, defined as 

$$ AG(n)=\{ x \mapsto a x +b \mid a,b \in \F_n,\, a\neq 0\}.$$

It is easy to check that the number of fixed points is in $\{0,1,n\}.$ Hence the nonzero $f_i$ for $i>1$ are $f_n$ and $f_{n-1}.$

A nongroup example of $2$-design can be obtained by considering

$$ \{ x \mapsto a x^3 +b \mid a,b \in \F_9,\, a\neq 0\}.$$

The conditions of the criterion can be checked in Magma \cite{M}. A computer free argument pointed out by a referee is that this set of permutations is in fact a one sided coset of
the affine linear group  $AG(n)$ by the Frobenius map $x\mapsto x^3.$ Such a coset is not a group ( as not containing the identity), and is a sharply transitive set of permutations, 
since $AG(n)$ is one.

A generalization of Theorem \ref{2crit} is as follows.

{\theorem \label{tcrit} A subset $D\subseteq S_n$ is a $t$-design in $(S_n,d_S),$ for $n\ge 2t,$ iff $$\sum_{i=1}^n f_i (\widehat{C_k(0)}-\widehat{C_k(i)})=\widehat{C_k(0)} $$ for $k \in [1..t].$ }

\begin{proof}
 The result follows by combining the vanishing of dual frequencies ( Proposition \ref{crit} ) with Lemma \ref{Charlier}.
\end{proof}

We are now ready for the main result of this paper.

{\theorem \label{sm} If $D$ is a $t$-design in $(S_n,d_S)$, then $|D| \geq n(n-1)\dots(n-t+1)$. In case of equality $f_i=0$ for $i\in [1..n-t].$ }

\begin{proof}
\begin{enumerate}
 \item By Burnside's Lemma applied to the action of $S_n$ on $t$-tuples
 \[\sum_{i=0}^{n}i(i-1)\dots(i-t+1)w_i = n!,\]
 or, equivalently
 
 \[\frac{1}{n!}\sum_{j=0}^{n}(n-j)(n-1-j)\dots(n-t+1-j)v_j =\langle 1, P_t\rangle_{S_n}=1,\]
 with $P_t(x)=(n-x)(n-1-x)\dots(n-t+1-x).$ By the $t$-design property, since $\deg(P_t)=t,$ we obtain $\langle 1, P_t\rangle_{S_n}=\langle 1, P_t\rangle_{D}.$
 Further, since $P_t(j)\ge 0,$ for $j \in [1..N]$ we see that $1=\langle 1, P_t\rangle_{D}\ge f_0 P_t(0).$ We conclude by  $f_0 = 1/|D|.$

 \item Apply Theorem 14.5.3 of \cite{G} with $\alpha=0$ and $p=P_t.$ The positivity of $p$ and the computation of $\langle 1, p\rangle$ are as above.
\end{enumerate}

\end{proof}

{\bf Remarks:} 
\begin{enumerate}
 \item This bound is best possible whenever sharply $t$-transitive sets of permutations on $n$ letters exist. Thus, in particular when $t=2$ and $n$ is a prime power.

 \item When $t=2$ this bound is strictly better than Corollary \ref{cor2}, a fact which shows that {\em tight} $2$-designs in the sense of \cite{G,Gp} do not exist. To define tight designs
 in general is beyond the scope of this paper; suffices to say that a tight 2-design would be one meeting Corollary 2 with equality.
 
 \item In case of equality the $t$ nonzero frequencies can be computed exactly by solving a linear system of order $t$ using Theorem \ref{tcrit}.
\end{enumerate}

%%%%%%%%%%%%%%%%%%%%%%%%%%%%%%%%%%%%%%%%%%%%%%%%%%%%%%%%%%%%%%%%%%%%%%%%%%%%%%%%%%%%%%%%%%%%%%%%%%%%%%%%%%%%%%%%%%%%%%%%%
\section{Conclusion}
%%%%%%%%%%%%%%%%%%%%%%%%%%%%%%%%%%%%%%%%%%%%%%%%%%%%%%%%%%%%%

In this article, we have given some bounds on the size of designs in the symmetric group. While $t$-transitive sets of permutations provide examples of such designs, it is not clear that
 all $t$-designs are of that form. We do not know any example of $t$-designs which are not $t$-transitive sets of permutations. 
 
 Another difficult open problem is the improvement of Theorem \ref{sm} for some specific values of $n$ and $t.$ For instance, showing that $|D|>90$ for a $2$-design in $S_{10}$ would give an alternative proof of
 the non-existence of the projective plane of order $10$ \cite{L,L+}. It is conceivable that adding linear constraints on the frequencies in the linear programming of Theorem 11 would lead to
 such an important result.
 
 %As for a general Delsarte MacWilliams bound for all DDR spaces, it seems that the theory is not ready yet for such a result.

{\bf Acknowledgement:} The authors thank  Andries Brouwer, Peter Cameron, John Cannon, Ferdinand Ihringer, Tatsuro Ito, Bill Martin and Sam Mattheus for helpful discussions.

\end{document}